\documentclass[preprint,11pt]{elsarticle}

\usepackage{amssymb}
\usepackage{multirow}
\usepackage{amssymb}
\usepackage{amsthm}
\usepackage{amsfonts}
\usepackage{amsmath}
\usepackage{graphicx}
\usepackage[top=16mm, bottom=16mm, left=16mm, right=16mm]{geometry}
\newcommand{\pq}[2]{\genfrac{[}{]}{0pt}{}{#1}{#2}}
\newcommand{\wt}{\mbox{wt}}

\newtheorem{theorem}{Theorem}
\newtheorem{lemma}[theorem]{Lemma}
\newtheorem{corollary}[theorem]{Corollary}
\newtheorem{proposition}[theorem]{Proposition}
\newdefinition{remark}{Remark}

\begin{document}

\begin{frontmatter}



\title{Rook theoretic proofs of some identities related to Spivey's Bell number formula}

\author[cnu]{R.B.~Corcino\corref{cor1}}
\ead{rcorcino@yahoo.com}
\author[upd]{R.O.~Celeste}
\ead{ching@math.upd.edu.ph}
\author[upd]{K.J.M.~Gonzales\fnref{fn1}}
\ead{kmgonzales@upd.edu.ph}

\cortext[cor]{Corresponding author}
\fntext[fn1]{Author is supported by Department of Science and Technology-Science Education Institute under its ASTHRD Program.}
\address[cnu]{Mathematics and ICT Department, Cebu Normal University, 6000 Cebu City, Philippines}
\address[upd]{Institute of Mathematics, University of the Philippines Diliman, 1101 Quezon City, Philippines}

\begin{abstract}
We use rook placements to prove Spivey's Bell number formula and other identities related to it, in particular, some convolution identities involving Stirling numbers and relations involving Bell numbers. To cover as many special cases as possible, we work on the generalized Stirling numbers that arise from the rook model of Goldman and Haglund. An alternative combinatorial interpretation for the Type II generalized $q$-Stirling numbers of Remmel and Wachs is also introduced in which the method used to obtain the earlier identities can be adapted easily.
\end{abstract}


\end{frontmatter}

\section{Introduction}\label{1}
Let $B(n)$ denote the $n$-th Bell number and $S(n,k)$ denote the Stirling number of the second kind. Spivey \cite{Spi} established combinatorially the following identity for Bell numbers
\begin{equation}\label{spv}
B(n+m) = \sum_{k=0}^n \sum_{j=0}^m j^{n-k} \binom{n}{k} S(m,j) B(k)\,.
\end{equation}
Various alternative proofs and extensions of this identity have appeared in the literature. For instance, Gould and Quaintance \cite{GouQua} provided a generating function proof, which, in turn, was extended by Xu \cite{Xu} in the case of Hsu and Shuie's \cite{HsuShu} generalized Stirling numbers. Katriel \cite{Kat} proved a $q$-analogue using certain $q$-differential operators. Still, Belbachir and Mihoubi \cite{BelMih} proved \ref{spv} using a decomposition of the Bell polynomial using a certain polynomial basis. See also \cite[Theorem 10]{MaaMih}, \cite[Theorem 5.3]{ManSchSha} and \cite{Mez} for other generalizations and methods. The present authors also derived a generalization \cite[Theorem 4.5]{CorCelGon} of Spivey's identity in the case of the generalized $q$-Stirling numbers that arise from normal ordering.

In this short paper, we derive variations of \ref{spv} and other related identities using the aforementioned rook model which we describe in Section \ref{2}. The results are presented in Section \ref{3}. In Section \ref{4}, we introduce a new rook model for the Type II generalized $q$-Stirling numbers of Remmel and Wachs \cite{RemWac} which is a variation of Goldman and Haglund's. We show how the method used in deriving the earlier identities can be easily adapted under this model.

\section{The rook model}\label{2}
In this section, we introduce a rook model based on Goldman and Haglund's \cite{GolHag}. We will later show that these two models are essentially identical and give rise to the same rook numbers.

\begin{figure}[htbp]
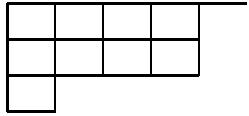

\begin{center}
\begin{tabular}{|c|c|c|c|c}
\cline{1-5}
\phantom{h} & \phantom{h} & \phantom{h} & \phantom{h} & \phantom{h} \\

\cline{1-4}
\phantom{h} &  &  &  \\

\cline{1-4} \phantom{h} \\

\cline{1-1}

\end{tabular}
\caption{The board outlined by $UVUUUVVU$.}\label{board}
\end{center}
\end{figure}

Let $w$ be a word consisting of the letters $U$ and $V$. If we let $U$ correspond to a unit horizontal step and $V$ a unit  vertical step, then $w$ outlines a Ferrers boards (or simply, board), which we denote by $B(w)$. For example, Figure \ref{board} shows the board outlined by $UVUUUVVU$. Note that we allow the rightmost columns and the bottommost rows to be empty, i.e., not containing any cell. The initial $U$'s and final $V$'s, if any, are extraneous when we consider rook placements, but they are natural in the context of normal ordering (see \cite{Var}). Alternatively, we can also describe a board by the length of its columns.

Given a board $B$, we associate to each cell of $B$ an integer called its \emph{pre-weight}. Fix $s\in\mathbb R$. A placement of $k$ rooks on a board $B$ is a marking of $k$ cells of $B$ with ``$\bullet$'' such that at most one rook is placed in each column. We assume that the rooks are placed from right to left among $k$ chosen column. Once a rook is placed on a cell, $s-1$ is added to the pre-weight of every cell to its left in the same row. However, if a cell lies above a rook (in which case it said to be \emph{cancelled} by the rook), then the cell is assigned the pre-weight 0. We denote by $C_k(B;s)$ the set of all placement of $k$ rooks on $B$ under the rules we just described.

For the purpose of distinguishing the pre-weights of a board before and after rooks are placed, we will call the former as \emph{default pre-weights}. We number the rows and columns of a board $B$ from top to bottom and from right to left, respectively.  The term pre-weight is chosen because it determines the weight of a board. Specifically, a cell with pre-weight $t$ is assigned the weight $q^t$ if it does not contain a rook, and $[t]_q=\frac{q^t-1}{q-1}$ if it contains a rook. The weight of a rook placement $\phi$, denoted by $\wt(\phi)$, is defined as the product of the weight of the cells.

\begin{figure}[htbp]
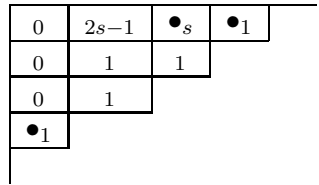

\begin{center}
\begin{tabular}{|c|c|c|c|c}
\cline{1-5}
${}_{0}$ & ${}_{2s-1}$ & $\bullet_s$ & $\bullet_1$ & \phantom{h} \\

\cline{1-4}
${}_0$ & ${}_1$ & ${}_1$  &  \multicolumn{1}{c}{\hphantom{x}} \\

\cline{1-3}
${}_0$ & ${}_1$ &   \multicolumn{1}{c}{\hphantom{x}} \\

\cline{1-2}
$\bullet_1$ &  \multicolumn{1}{c}{\hphantom{x}} \\

\cline{1-1} \multicolumn{1}{|c}{\hphantom{x}}  \\

\end{tabular}
\caption{A rook placement on $J_5$.}\label{rookplc}
\end{center}
\end{figure}

Denote by $J_n$ the board outlined by $(VU)^n$ whose cells have default pre-weight 1. Note that columns of $J_n$ have lengths $0, 1, \ldots n-1$. Figure \ref{rookplc} shows a rook placement on $J_5$ where the values in the cells indicate the pre-weights. This particular rook placement has weight $q^{2s+2}[s]_q$.

We denote by $S_{s,q}[n,k]$ the sum of the weights of all placements of $n-k$ rooks on $J_n$, i.e.,
\[
S_{s,q}[n,k] = \sum_{\phi\in C_{n-k}(J_n;s)} \wt(\phi)\,.
\]

\begin{remark} If, for instance, a column has two cells where the bottom cell has pre-weight $c_1$ and the top cell has pre-weight $c_2$, then the placement of a rook on the bottom cell has weight $[c_1]_q$ while the placement of a rook on the top cell has weight $q^{c_1}[c_2]_q$. Hence, the total weight of all rook placements on this column is $[c_1]_q+q^{c_1}[c_2]_q=[c_1+c_2]_q$. In general, if $c$ is the sum of the pre-weights of the cells in a column, then the total weight of all rook placements on the column is $[c]_q$. \label{ddd}
\end{remark}

\begin{proposition}\label{trr} The number $S_{s,q}[n,k]$ satisfies the recurrence
\begin{equation*}
S_{s,q}[n,k] = q^{s(n-1)-(s-1)(k-1)} S_{s,q}[n-1,k-1] + [s(n-1)-(s-1)k]_q S_{s,q}[n-1,k]\,,
\end{equation*}
with initial conditions $S_{s,q}[n,0]=S_{s,q}[0,n]=\delta_{0,n}$, where $\delta_{0,n}=1$ if $n=0$ and 0 if $n\neq 0$.
\begin{proof}If there is a rook on the $n$-th column, then the other columns form a rook placement from $C_{n-1-k}(J_{n-1};s)$. Due to the placement of $n-1-k$ rooks on the first $n-1$ columns, the $n$-th column has pre-weight $(n-1)+(s-1)(n-1-k)=s(n-1)-(s-1)k$. Hence, the total weight contributed by all rook placements on the $n$-th column is $[s(n-1)-(s-1)k]_q$.

If there is no rook on the $n$-th column, then the remaining columns form a rook placement from $C_{n-k}(J_{n-1};s)$. Because of the placement of $n-k$ rooks on the first $n-1$ columns, the $n$-th column has pre-weight $(n-1)+(s-1)(n-k)=s(n-1)-(s-1)(k-1)$, which implies that the $n$-th column contributes a weight of $q^{s(n-1)-(s-1)(k-1)}$.
\end{proof}
\end{proposition}

By comparing recurrences, we see that $h^{n-k}S_{s,q}[n,k]$ equals the number $\mathfrak S_{s;h}(n,k|q)$ in \cite{ManSchSha2}. These numbers are coefficients of the string $(VU)^n$ in its normally ordered form, i.e., an equivalent expression where the $V$'s are to the right of $U$'s obtained using the relation $UV=qVU+hV^s$. When $h=1,q=1$, the cases $s=1$ and $s=0$ produce the usual Stirling number of the first kind and Stirling number of the second kind, respectively. (See \cite{ManSchSha2} for other special cases and their relationship with normal ordering.)

\begin{remark} The original Goldman-Haglund rook model is as follows (see \cite[Section 7]{GolHag}). Let $B$ be a board and consider a rook placement $\phi$ on $B$ such that no two rooks occupy the same column. Given a cell $\gamma$ in $B$, let $v(\gamma)$ be the number of rooks strictly to the right of, and in same row as $\gamma$. Define the weight of $\gamma$ to be

\begin{align*}
\wt(\gamma) = \begin{cases}
 1 & \mbox{if there is a rook below} \\ & \hspace{0.3in}\mbox{and in the same column as $\gamma$,} \\
 [(s-1)v(\gamma)+1]_q & \mbox{if $\gamma$ contains a rook,}\\
q^{(s-1)v(\gamma)+1}& \mbox{else}\,.
\end{cases}
\end{align*}

The weight of $\phi$, denoted by $\wt(\phi)$ is then defined as the product of the weight of the cells. (Note that in \cite{GolHag}, the boards are bottom and right justified, the rooks are placed from left to right and a rook cancels cells that lie below it. In addition, $\alpha$ is used instead of $s$. The definition of $\wt(\gamma)$ above has been modified to reflect the definitions in the current paper.)

It can be easily seen that the weight of a rook placement under the Goldman-Haglund model and under the current model are identical. First, note that $\wt(\gamma)$ is precisely the pre-weight of $\gamma$. Under the latter, a cell lying above a rook has pre-weight 0 and hence, has weight $q^0=1$. A cell containing a rook and has $t$ rooks lying to its right and on the same row will have pre-weight $t(s-1)+1$ and hence, the weight $[t(s-1)+1]_q$. Finally, a cell not containing a rook and has $t$ rooks lying to its right and on the same row will have pre-weight $t(s-1)+1$ and thus, the weight $q^{t(s-1)}+1$. Hence, given a rook placement $\phi$, $\wt(\phi)$ is identical under both models.
\end{remark}

The concept of pre-weights provides us with some degree of flexibility in further generalizing the earlier model. For instance, instead of adding $s-1$ to the pre-weight of each cell to the left of a rook, we can fix a \emph{rule} $R$ which specifies, given a placement of rook in a cell, the cell (in each column to the right of a rook) that receives the pre-weight increment. Denote by $C^R_k(B;s)$ the set of placement of $k$ rooks on the board $B$ under this rule. In particular, we denote by $\mathcal R$ the rule which specifies that the $i$-th rook (from the right) adds a pre-weight of $s-1$ to the $i$-th cell above the bottom cell in every column to its left. Of course, when the default rule is used, we simply use the notation $C_k(B;s)$.

\begin{figure}[htbp]
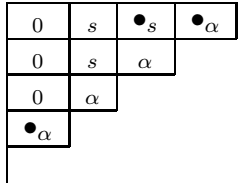

\begin{center}
\begin{tabular}{|c|c|c|c|c|}
\cline{1-4}
${}_{0}$ & ${}_{s}$ & $\bullet_{s}$ & $\bullet_{\alpha}$ \\

\cline{1-4}
${}_0$ & ${}_{s}$ & ${}_{\alpha}$  &  \multicolumn{1}{c}{\hphantom{x}} \\

\cline{1-3}
${}_0$ & ${}_{\alpha}$ &   \multicolumn{1}{c}{\hphantom{x}} \\

\cline{1-2}
$\bullet_{\alpha}$ &  \multicolumn{1}{c}{\hphantom{x}} \\

\cline{1-1} \multicolumn{1}{|c}{\hphantom{x}}  \\

\end{tabular}
\caption{A rook placement on $J'_{4,\alpha}$ using rule $\mathcal R$.}\label{rookplc2}
\end{center}
\end{figure}

In addition to specifying a rule $R$, we can also assign other values for the default pre-weights. This will be useful when we consider Type II $q$-analogues of Stirling numbers in Section \ref{4}. For $\alpha\in\mathbb R$, denote by $J'_{n,\alpha}$ the board outlined by $(VU)^nV$ where the bottom cell in each column each has default pre-weight $\alpha$. Figure \ref{rookplc2} shows one rook placement in $C^{\mathcal R}_{3}(J'_{4,\alpha};s)$, which has weight $q^{2s+2\alpha}[\alpha]_q^2[s]_q$.

\begin{proposition}\label{prew} Given any rule $R$ and two boards $B_1$ and $B_2$ with the same number of non-empty columns such that the corresponding columns have identical total column pre-weights, we have $\sum_{\phi\in C^{R}_k(B_1;s)}  \wt(\phi) = \sum_{\phi\in C^{R}_k(B_2;s)}  \wt(\phi)$. Furthermore, given a board $B$ and two rules $R_1$ and $R_2$, $\sum_{\phi\in C^{R_1}_k(B;s)}  \wt(\phi) = \sum_{\phi\in C^{R_2}_k(B;s)}  \wt(\phi)$.
\begin{proof} We can establish, using the same reasoning as Remark \ref{ddd}, that the sum of the weights of all possible rook placements on a particular column is completely determined by the column's total pre-weight, regardless of the distribution of the pre-weights of the individual cells. This proves the first statement. The second statement follows from the first statement and the fact that different rules applied on the same board result to boards with identical total column pre-weights.
\end{proof}
\end{proposition}

\section{Results}\label{3}
Define the $n$-th generalized Bell polynomial by \[B_{s,q}[n;x]=\sum_{k=0}^n S_{s,q}[n,k] x^k\,.\] The $n$-th generalized Bell number, denoted by $B_{s,q}[n]$, is defined as $B_{s,q}[n;1]$. The $q$-binomial coefficients are given by $\pq{n}{k}_q=\frac{[n]_q!}{[k]_q![n-k]_q!}$, where $[n]_q!=[n]_q[n-1]_q\cdots[2]_q[1]_q$. Clearly, $\pq{n}{k}_q = \pq{n}{n-k}_q$. They also satisfy the relations (see \cite[Table 1 and Identity (2.2)]{MedLer})
\begin{align*}
\pq{n}{k}_q &= \sum_{0\leq t_1 \leq t_2 \leq \ldots \leq t_{n-k}\leq k} q^{t_1+t_2+\cdots+t_{n-k}}\\
\pq{n}{k}_q &= \sum_{t_0+t_1+t_2+\cdots+t_{k}=n-k} q^{0t_0+1t_1+2t_2+\cdots+kt_{k}}\,.
\end{align*}

In this section, we show how Spivey's identity, specifically, its generalization involving $S_{s,q}[n,k]$, is a consequence of a convolution identity. We then proceed by deriving other convolution identities and the identities for Bell numbers that arise from them.

A version of Lemma \ref{lemma1} and Theorem \ref{or} already appeared in \cite{CorCelGon}, but only the case where $s\in\mathbb N$ and with pre-weights interpreted as subdivisions in cells was proved in detail.

\begin{lemma}\label{lemma1} Let $\phi$ be a placement of $k$ rooks on the board $J_n$. Then, there exists a unique (possibly empty) collection $\mathcal C$ of columns in $\phi$, such that if $|\mathcal C|=\mu+1$, then (a) these columns have a rook in the bottom $1,2,\ldots,\mu+1$ cells and (b) every column not in $\mathcal C$ contains at least $1+t$ uncanceled cells not containing a rook, where $t$ is the number of columns in $\mathcal C$ to the right of that column.
\begin{proof}
Let $\phi\in C^{R}_k(J'_{n,\alpha};s)$. The elements of $\mathcal C$ can be obtained iteratively as follows. Let $c_1$ be the first column from the right containing a rook on the bottom cell. If there exists no such $c_1$, then $\mathcal C=\varnothing$. Let $c_2$ be the next column containing a rook in one of the bottom two cells, i.e., either on the bottom cell or on the cell above the bottom cell. If there exists no such $c_2$, then $\mathcal C=\{c_1\}$. We can continue this process as long as needed until the elements of $\mathcal C$ are all determined. Note that this process also shows the uniqueness of $\mathcal C$.
\end{proof}
\end{lemma}

\begin{lemma}\label{alpha} Let $\alpha\in\mathbb R,n,k\in\mathbb N$. Then, for any rule $R$,
\begin{align}
\sum_{\phi\in C^R_{n-k}(J'_{n;\alpha};s)} \wt(\phi) &= \sum_{r=0}^n q^{\alpha r}\pq{n}{r}_{q^s} S_{s,q}[r,k] \prod_{i=0}^{n-r-1}[\alpha+si]_q \label{oh}\\
\sum_{\phi\in C^R_{n-k}(J'_{n;\alpha};s)} \wt(\phi) &= \sum_{r=0}^n q^{\alpha k} S_{s,q}[n,r] \pq{r}{k}_{q^{s-1}} \prod_{i=0}^{r-k-1} [\alpha+i(s-1)]_q\,.\label{oh1}
\end{align}
\begin{proof}
We prove \ref{oh} first. Let $\phi\in C^{\mathcal R}_{n-k}(J'_{n;\alpha};s)$. By Lemma \ref{lemma1}, we can write $\phi$ uniquely as a pair $(\mathcal C_{\phi},\phi-\mathcal C_{\phi})$, where $\mathcal C_{\phi}$ are the cells of $\phi$ that either lie in the columns that satisfy the properties in the lemma with $\mu=n-r-1$ or on the bottom $1+t$ cells of the other columns, and $\phi-\mathcal C_{\phi}$ are the cells in $\phi$ not in $\mathcal C_{\phi}$.  One sees that $\phi-\mathcal C_{\phi}$ forms some rook placement in $C_{r-k}(J_r;s)$. This accounts for the factor $S_{s,q}[r,k]$.

Let $L_{\mu}$ be the set of all such possible $\mathcal C_{\phi}$, with $\phi\in C^{\mathcal R}_{n-k}(J'_{n;\alpha};s)$ with $\mu=n-r-1$ as in Lemma \ref{lemma1}. It suffices to compute the sum of the weights of the elements of $L_{\mu}$. Clearly, the columns containing rooks contribute a weight of $\prod_{i=0}^{n-r-1} [\alpha+si]_q$ and the bottom $r$ cells of the columns not containing rooks contribute a weight of $q^{\alpha r}$ regardless of where the rooks are. On the other hand, the weight contributed by the remaining cells depends on the placement of the rooks, and varies as $q^{st_1} q^{st_2} \cdots q^{st_{n-r}}$ for some $0\leq t_1\leq t_2 \leq \ldots \leq t_{n-r} \leq r$. Hence, these cells contribute
\[
\sum_{0\leq t_1\leq t_2 \leq \ldots \leq t_{n-r} \leq r} q^{st_1} q^{st_2} \cdots q^{st_{n-r}} = \pq{n}{r}_{q^s}\,.
\]

For \ref{oh1}, we consider rook placements using the default rule on the board outlined by $(VU)^nV$ where the cells in the first row each have pre-weight $\alpha$ while all the other cells have pre-weight 1. By Proposition \ref{prew}, the sum of the weight of all placements of $n-k$ rook on this board equals that of $C^R_{n-k}(J'_{n;\alpha};s)$ for any rule $R$.

A rook placement in $C_{n-k}(J'_{n;\alpha};s)$ can then be formed as follows. First, we place $n-r$ rooks in cells lying on rows $2,\ldots,n-1$ and columns $2,\ldots,n$. The sum of the weights of all such placements is $S_{s,q}[n,r]$. Then, place the remaining $r-k$ rooks on the first row. The weight contributed by the cells in the first row containing rooks is $\prod_{i=0}^{r-k-1} [\alpha+i(s-1)]_q$, and this is independent of the choice of columns. On the other hand, the weight contributed by the cells in the first row not containing rooks is $q^{\alpha t_0} q^{\alpha t_1}\cdots q^{\alpha t_{r-k}}$ for some $t_0+t_1+\ldots+t_{r-k}=k$. Finally,
\begin{align*}
\sum_{t_0+t_1+\ldots+t_{r-k}=k} &q^{\alpha t_0} q^{\alpha t_1}\cdots q^{\alpha t_{r-k}} = q^{\alpha } \pq{r}{r-k}_{q^{s-1}} = q^{\alpha } \pq{r}{k}_{q^{s-1}}\,.
\end{align*}
\end{proof}
\end{lemma}

\begin{remark}
It is interesting to note that the individual terms on the $RHS$ of \ref{oh} and \ref{oh1} are not equal for fixed $r$.
\end{remark}

\begin{theorem}\label{or} Let $n,m,k\in\mathbb N$. Then,
\begin{align}\label{iden1}
S_{s,q}[n+m,k] = \sum_{r=0}^n \sum_{j=0}^m S_{s,q}[m,j]\,q^{r(j(1-s)+sm)} \pq{n}{r}_{q^s} S_{s,q}[r,k-j] \prod_{i=0}^{n-r-1} [j(1-s)+sm+si]_q\,.
\end{align}
Consequently,
\begin{equation}\label{bbb1}
B_{s,q}[n+m;x] = \sum_{r=0}^n \sum_{j=0}^m S_{s,q}[m,j]\,q^{r(j(1-s)+sm)} \pq{n}{r}_{q^s} B_{s,q}[r;x] x^j \prod_{i=0}^{n-r-1} [j(1-s)+sm+si]_q\,.
\end{equation}
In particular,
\begin{align}\label{an}
B_{s,q}[n+m] &= \sum_{r=0}^n \sum_{j=0}^m S_{s,q}[m,j]\,q^{r(j(1-s)+sm)} \pq{n}{r}_{q^s} B_{s,q}[r] \prod_{i=0}^{n-r-1} [j(1-s)+sm+si]_q\,.
\end{align}

\begin{proof} The number $S_{s,q}[n+m,k]$ equals the sum of the weights of all rook placements in the set $C_{n+m-k}(J_{n+m;s})$. The rooks may be placed as follows. First, place $m-j$ rooks in columns $2,\ldots,m$. The sum of the weights of all such placements is $S_{s,q}[m,j]$. Next, place the remaining $n+j-k$ rooks in columns $m+1,\ldots,n$. Due to the placement of the first $m-j$ rooks, columns $m+1,\ldots,n$ form a board whose column pre-weights are identical to those of $J'_{n;j(1-s)+sm}$. Identity \ref{iden1} follows from \ref{oh} of Theorem \ref{alpha} using the replacement $\alpha \rightarrow j(1-s)+sm, k \rightarrow k-j$.

To obtain \ref{bbb1}, multiply both sides of \ref{iden1} by $x^k$ and take the sum over all $0\leq k \leq n+m$.
\end{proof}
\end{theorem}


\begin{theorem}\label{ne} Let $n,m\in\mathbb N$. Then,
\begin{align}
S_{s,q}[n+m,k] = \sum_{r=0}^n \sum_{j=0}^m S_{s,q}[m,j] q^{(j(1-s)+sm)(k-j)} \pq{r}{k-j}_{q^{s-1}} S_{s,q}[n,r] \prod_{i=0}^{r-k+j-1} [j(1-s)+sm+i(s-1)]_q\,.
\end{align}
Consequently,
\begin{align}
B_{s,q}[n+m;x] &=\sum_{r,j,k} S_{s,q}[m,j] q^{(j(1-s)+sm)(k-j)} \pq{r}{k-j}_{q^{s-1}} S_{s,q}[n,r] x^{k} \prod_{i=0}^{r-k+j-1} [j(1-s)+sm+i(s-1)]_q\,.
\end{align}
In particular,
\begin{align}\label{exd}
B_{s,q}[n+m] &=\sum_{r,j,k} S_{s,q}[m,j] q^{(j(1-s)+sm)(k-j)} \pq{r}{k-j}_{q^{s-1}} S_{s,q}[n,r] \prod_{i=0}^{r-k+j-1} [j(1-s)+sm+i(s-1)]_q\,.
\end{align}
\begin{proof} The proof is the same as that of Theorem except that we use Identity \ref{oh1} of Theorem \ref{alpha}.
\end{proof}
\end{theorem}

Taking cue from the identities in \cite[Theorem 2.6]{MedLer}, we derive another set of convolution identities.

\begin{theorem}\label{sec}For $n,m,j\in\mathbb N$, we have
\begin{align}
S_{s,q}[n+1,m+j+1] &= \sum_{k=m}^n \sum_{r=0}^{n-k} S_{s,q}[k,m]\, q^{r((k-m)(s-1)+k+1)+k+(k-m)(s-1)} \nonumber \\ &\hspace{0.5in} \pq{n-k}{r}_{q^s} S_{s,q}[r,j] \prod_{i=0}^{n-k-r-1} [(k-m)(s-1)+k+1+si]_q\,,\label{hey1} \\
S_{s,q}[n+1,m+j+1] &= \sum_{k=m}^n \sum_{r=0}^{n-k} S_{s,q}[k,m]\, q^{(j + 1)(k + (k - m)(s - 1)) + j} \nonumber \\ &\hspace{0.5in} \pq{r}{j}_{q^{s-1}} S_{s,q}[n-k,r] \prod_{t=0}^{r-j-1} [(k-m)(s-1)+k+1+t(s-1)]_q\,.\label{hey2}
\end{align}
\begin{proof} For a given rook placement in $C_{n-m-j}(J_{n+1};s)$, there exists a unique $k$ such that there are $k-m$ rooks in columns $2,\ldots,k$, $0$ rooks in column $k+1$ and $n-j-k$ rooks in columns $k+2,k+3,\ldots,n+1$. Indeed, if there exists $k_1<k_2$ satisfying these properties, then there are $k_2-k_1$ rooks in columns $k_1+1,k_1+2,\ldots,k_2$. This is impossible since there is no rook in column $k_1+1$.

Suppose that $k-m$ rooks have been placed in columns $1,2,\ldots,k$. The weight contributed by the $k+1$-th column is $q^{k+(k-m)(s-1)}$. Using Theorem \ref{alpha} with the substitutions $\alpha\rightarrow k+1+(k-m)(s-1), n\rightarrow n-k, k\rightarrow j, i\rightarrow r$, we obtain two different expressions for the sum of the weights of the placements of $n-j-k$ rooks in columns $k+2,k+3,\ldots,n+1$. 
\end{proof}
\end{theorem}


\begin{remark}We briefly discuss some special cases.

By letting $s=0$ in \ref{an}, we get the following $q$-analogue of \ref{spv} which was previously derived by Katriel \cite{Kat}
\[
B_q[n+m]=\sum_{r=0}^n \sum_{j=0}^m S_q[m,j] q^{rj} \binom{n}{r} B_q[r] [j]_q^{n-r}\,,
\]
where, $S_{0,q}[n,k]=S_q[n,k]$ and $B_{0,q}[n+m]=B_{q}[n+m]$ are the $q$-Stirling number of the second kind and the $q$-Bell number, respectively.

A ``dual'' of Spivey's identity was derived by Mez\"o \cite{Mez} and is given by
\begin{align}\label{mez}
(n+m)! = \sum_{r=0}^n \sum_{j=0}^m c(m,j) m^{\overline{n-r}} \binom{m}{j} r!\,,
\end{align}
where $a^{\overline{b}}=a(a+1)\cdots(a+b-1)$ and $c(n,k)$ is the Stirling number of the first kind. One can check that\ref{an} reduces to \ref{mez} when $s=1,q=1$. Using \ref{exd}, we deduce another expression for $(n+m)!$, namely
\[
(n+m)! = \sum_{r,j,k} c(m,j) \binom{r}{k-j} c(n,r) m^{r-k+j}\,,
\]
or equivalently,
\[
(n+m)! = \sum_{r,j,k} c(m,j) \binom{r}{k} c(n,r) m^{k}\,.
\]
A simple combinatorial interpretation of the above identity was given by Mez\"o \cite{Mez2}. Color the first $m$ elements red and the remaining $n$ elements blue. Place the red elements into $j$ cycles in $c(m,j)$ ways and the blue elements into $r$ cycles in $c(n,r)$ ways. Then, write the cycles in standard form, i.e., such that the maximum element in each cycle is written first and the cycles are arranged by increasing first elements. Pick $k$ cycles among the $r$ blue cycles in $\binom{r}{k}$ ways. Insert the $k$ blue cycles among the red elements such that if there are consecutive blue cycles after a red element, we rewrite the blue cycles in standard form. This can be done in $m^k$ ways. Finally, we remove the parentheses of the blue cycles which were inserted and forget the colors.

Other identities for the the classical Stirling numbers and Bell numbers and their $q$-analogues may be obtained as straight-forward corollaries of the results in this section. For instance, Theorem \ref{ne} gives
\begin{align*}
S_q[n+m,k] &= \sum_{r=0}^n \sum_{j=0}^m S_q[m,j] q^{j(k-j)} \binom{r}{k-j}_{1/q} S_q[n,r] [j]_q[j-1]_q\cdots[k-r+1]_q\,,
\end{align*}
which appears to be new.
\end{remark}

The proposition that follows generalizes Theorems \ref{or}, \ref{ne} and \ref{sec}. As before, two explicit forms, which we no longer state, can be obtained through an application of Theorem \ref{alpha}. For the proofs, only the needed partition of the boards are described.

\begin{proposition}There holds 
\begin{align*} 
S_{s,q}[m_1+\cdots+m_n,k] &=\sum_{\substack{j_1+\cdots+j_n=k\\ j_1,\ldots,j_n\in\mathbb N}}\,\,
\prod_{i=1}^n \sum \wt(\phi)\,,\\
S_{s,q}[n+j-1,m_1+\cdots+m_j+j-1] &= \sum_{\substack{k_1+\cdots+k_j=n \\ k_1,\ldots,k_n\in\mathbb N}} q^{\sum_{t=1}^{i-1} k_t+(k_t-m_t)(s-1)} \prod_{i=1}^j \sum \wt(\phi) \,.
\end{align*}
where in the first expression, the rightmost sum is taken over all rook placements in $C_{m_i-j_i}\left(J'_{m_i;\sum_{t=1}^{i-1} j_t(1-s)+sm_t};s\right)$ and in the second expression, $C_{k_i-m_i}\left(J'_{k_i;\sum_{t=1}^{i-1} k_t+1+(k_t-m_t)(s-1)};s\right)$.
\begin{proof} For the first identity: starting from the right, divide the the board $J_{m_1+\cdots+m_n}$ into $n$ groups $G_1,\ldots,G_n$ consisting of $m_1,\ldots,m_n$ adjacent columns, respectively. Then, for $1\leq i\leq n$, place $m_i-j_i$ rooks on the columns in $G_i$.

For the second identity: use the fact that for every rook placement in $C_{n-(m_1+\cdots+m_j)}(J_{n+j-1})$, there exists a unique $j$-tuple $(k_1,\ldots,k_j)$ satisfying $k_1+\cdots+k_j=n$ such that (a) starting from the right, the board $J_{n+j-1}$ can be divided into $j$ groups $G_1,\ldots,G_j$ consisting of $k_1,\ldots,k_n$ adjacent columns, respectively, with one column not containing a rook in between, and (b)  for $1\leq i\leq j$, there are $k_i-m_i$ rooks on the columns in $G_i$.
\end{proof}
\end{proposition}


\section{Extension to Type II $q$-analogues}\label{4}
Hsu and Shuie \cite{HsuShu} introduced generalized Stirling numbers via the relations
\begin{align*}
(x|\alpha)_n &=\sum_{k=0}^n S^1_{n,k}(\alpha,\beta,\rho)(x-\rho|\beta)_k\\
(x|\beta)_n &=\sum_{k=0}^n S^2_{n,k}(\alpha,\beta,\rho)(x+\rho|\alpha)_k
\end{align*}
where $(z|\alpha)_0=1$ and $(z|\gamma)_n=z(z-\gamma)\cdots(z-(n-1)\gamma)$ for every positive integer $n$. When $(\alpha,\beta,r)=(0,1,0)$, $S^1_{n,k}(\alpha,\beta,\rho)$ and $S^2_{n,k}(\alpha,\beta,\rho)$ become the usual Stirling number of the first kind and Stirling number of the second kind, respectively. Since $S^1_{n,k}(\alpha,\beta,\rho)=S^2_{n,k}(\beta,\alpha,-\rho)$, it suffices to consider just one of these numbers, namely $S^1_{n,k}(\alpha,\beta,\rho)$.

Let $[\gamma]_{p,q}=\frac{p^{\gamma}-q^{\gamma}}{p-q}$ for any $\gamma\in\mathbb R$. The Type II $p,q$-analogues of these numbers were introduced by Remmel and Wachs \cite{RemWac} by the following recursions ($\rho$ appears as $-r$ in \cite{RemWac})
\begin{align*}
S^{1,p,q}_{n,k}(\alpha,\beta,\rho) &= q^{(k-1)\beta-(n-1)\alpha+\rho} S^{1,p,q}_{n-1,k-1}(\alpha,\beta,\rho) + p^{t-k\beta} [k\beta-(n-1)\alpha+\rho] S^{1,p,q}_{n-1,k}(\alpha,\beta,\rho)\\
S^{2,p,q}_{n,k}(\alpha,\beta,\rho) &= q^{-\rho+(k-1)\alpha-(n-1)\beta} S^{2,p,q}_{n-1,k-1}(\alpha,\beta,\rho) + p^{t+\rho-k\alpha} [k\alpha-\rho-(n-1)\beta] S^{2,p,q}_{n-1,k}(\alpha,\beta,\rho)
\end{align*}
where in addition $S^{1,p,q}_{0,0}(\alpha,\beta,\rho)=S^{2,p,q}_{0,0}(\alpha,\beta,\rho)=1$ and $S^{1,p,q}_{n,k}(\alpha,\beta,\rho)=S^{2,p,q}_{n,k}(\alpha,\beta,\rho)=0$ if $k=0$ or $k>n$. Since $S^{1,p,q}_{n,k}(\alpha,\beta,\rho)=S^{2,p,q}_{n,k}(\beta,\alpha,-\rho)$, we also consider just one of these numbers, which we choose to be $S^{1,p,q}_{n,k}(\alpha,\beta,\rho)$. When $(\alpha,\beta,\rho)=(0,1,0)$ and $p=q=1$, $S^{1,p,q}_{n,k}(\alpha,\beta,\rho)=S(n,k)$.

Remmel and Wachs also gave combinatorial interpretations for these numbers for a certain choice of parameters. Our goal in this section is to introduce an alternative rook interpretation for $S^{1,1,q}_{n,k}(\alpha,\beta,\rho)$ which adapts the method used in Section \ref{3}. To do this, we will derive the analogue of Theorem \ref{or} and Theorem \ref{ne}. We leave the corresponding versions of the other identities in Section \ref{3} to the interested reader.

Let $c,d\in\mathbb R$ and denote by $J^{c,d}_n$ the board outlined by $(VU)^n$ such that each bottom cell has default pre-weight $d$ and all other cells have default pre-weight $c+d$. Let $S^{c,d}_{s,q}[n,k]$ denote the sum of the weights of all placements of $n-k$ rooks on $J^{c,d}_n$.

\begin{remark} Since $J^{1,0}_{n}=J_n$, it follows that $S^{1,0}_{s,q}[n,k]=S_{s,q}[n,k]$.
\end{remark}

\begin{remark} A ``special case'' of $J^{c,d}_{n-1}$ with $c=m,d=0$ are the $m$-jump boards $Jb_{n,m}$ \cite[Example 3]{GolHag} which have column heights $0,m,2m,\ldots,(n-1)m$. One sees that these two boards have identical column pre-weights.
\end{remark}

\begin{proposition} The number $S^{c,d}_{s,q}[n,k]$ satisfies the recurrence
\[
S^{c,d}_{s,q}[n,k] = q^{c(n-1)+d+(n-k)(s-1)} S^{c,d}_{s,q}[n-1,k-1] + [c(n-1)+d+(n-k-1)(s-1)]_q S^{c,d}_{s,q}[n-1,k]
\]
Hence,
\[
S^{\beta-\alpha,-\rho}_{-\beta+1,q}[n,k] = S^{1,1,q}_{n,k}(\alpha,\beta,\rho)\,.
\]
\begin{proof} The proof is similar to that of Proposition \ref{trr}.\end{proof}
\end{proposition}

\begin{theorem}\label{t1}
For $n,m,k\in\mathbb N$, we have
\begin{multline*}
S^{c,d}_{s,q}[n+m,k] = \sum_{r=0}^n \sum_{j=0}^m S^{c,d}_{s,q}[m,j]\, q^{r(d+mc+(m-j)(s-1))} \\ \pq{n}{r}_{q^{c+s-1}} S^{c,0}_{s,q}[r,k-j] \prod_{i=0}^{n-r-1} [d+mc+(m-j)(s-1)+(c+s-1)i]_q\,.
\end{multline*}
This implies that
\begin{multline*}
S^{1,1,q}_{n,k}(\alpha,\beta,\rho) = \sum_{r=0}^n \sum_{j=0}^m S^{1,1,q}_{m,j}(\alpha,\beta,\rho)\, q^{r(\beta j-\alpha m-\rho)} \pq{n}{r}_{q^{-\alpha}} S^{1,1,q}_{r,k-j}(\alpha,\beta,0) \prod_{i=0}^{n-r-1} [\beta j-\rho-\alpha(m+i)]_q\,.
\end{multline*}
\begin{proof}
We proceed as in Theorem \ref{or}. The $m-j$ rooks in columns $2,\ldots,m$ explains the factor $S^{c,d}_{s,q}[m,j]$. We then determine the sum of the weights of all placements of $n+j-k$ rooks in the remaining columns.

The placement of $m-j$ rooks adds $(m-j)(s-1)$ to the total pre-weights of columns $m+1,\ldots,n$. These columns form a board whose column pre-weights are identical to the board outlined by $(VU)^nV$ such that the bottom cells have pre-weight $d+mc+(m-j)(s-1)$ and the other cells have pre-weight $c$. For convenience, let us denote such board by $B$.

Note that Lemma \ref{lemma1} is a statement regarding only the placement of rooks and thus, it holds true regardless of the default pre-weights of the cells and the rule used. Hence, we can use this lemma on $B$ with rule $\mathcal R$, reasoning as in Theorem \ref{or}.

Let $\phi\in C^{\mathcal R}_{n+j-k}(B;s)$. Again, by Lemma \ref{lemma1}, we can write $\phi$ uniquely as $(\mathcal C_{\phi},\phi-\mathcal C_{\phi})$. Note that $\phi-\mathcal C_{\phi}$ forms a rook placement in $C_{r-k+j}(J^{c,0}_r;s)$ and not in $C_{r-k+j}(J^{c,d}_r;s)$! This is because $\phi-\mathcal C_{\phi}$ does not include any bottom cell of $\phi$.

The board $B$ is similar to $J'_{n;\alpha}$ with $\alpha=d+mc+(m-j)(s-1)$, the only difference is that the non-bottom cells of $B$ have pre-weight $c$. Hence, a cell with pre-weight $c$ which receives an additional pre-weight due to the placement of a rook following rule $\mathcal R$ will have pre-weight $c+s-1$.
\end{proof}
\end{theorem}

Define $B^{1,1,q}_{n;x}(\alpha,\beta,\rho)=\sum_{k=0}^n S^{1,1,q}_{n}(\alpha,\beta,\rho)$.
\begin{corollary}
For $n,m,\in\mathbb N$, we have
\begin{multline}\label{mezz}
B^{1,1,q}_{n+m;x}(\alpha,\beta,\rho) = \sum_{r=0}^n \sum_{j=0}^m S^{1,1,q}_{m,j}(\alpha,\beta,\rho)\, q^{r(\beta j-\alpha m-\rho)} \pq{n}{r}_{q^{-\alpha}} B^{1,1,q}_{r;x}(\alpha,\beta,0) \prod_{i=0}^{n-r-1} [\beta j-\rho -\alpha(m+i) ]_q\,.
\end{multline}
\end{corollary}

Letting $q=1$ produces a generalization of Spivey's identity, which is a variant of Xu's result \cite[Corollary 8]{Xu}. Also, Mez\"o \cite[Theorem 2]{Mez} obtained a special case of \ref{mezz} for $(\alpha,\beta,\rho)=(0,1,\rho)$.

\begin{theorem} Let $n,m\in\mathbb N$. Then,
\begin{multline*}
S^{c,d}_{s,q}[n+m,k] = \sum_{r=0}^n \sum_{j=0}^m S^{c,d}_{s,q}[m,j] q^{(d + mc + (m - j)(s - 1))(k - j)} \\ \pq{r}{k-j}_{q^{s-1}} S^{c,0}_{s,q}[n,r] \prod_{i=0}^{r-k+j-1} [d + mc + (m - j)(s - 1) + i(s-1)]_q\,.
\end{multline*}
This gives
\begin{multline}
S^{1,1,q}_{n+m,k}(\alpha,\beta,\rho) = \sum_{r=0}^n \sum_{j=0}^m S^{1,1,q}_{m,j}(\alpha,\beta,\rho) q^{(\rho+\alpha m -\beta j)(j-k)} \\ \pq{r}{k-j}_{q^{-\beta}} S^{1,1,q}_{n,r}(\alpha,\beta,0) \prod_{i=0}^{r-k+j-1} [\rho+\alpha m -\beta (j+i)]_q\,.
\end{multline}
\begin{proof} Instead of using the board $J^{c,d}_{n}$, we use the board outlined by $(VU)^n$ such that the bottom cells of columns $2,\ldots,m$ have pre-weight $c+d$, the cells in row $m$ have pre-weight $c+d$, and all other cells have pre-weight $c$. We then proceed as in the proof of Theorem \ref{ne}.
\end{proof}
\end{theorem}
\section{Some Remarks}
Mez\"o's dual of Spivey's identity is based on the relation $n!=\sum_{k=0}^n c(n,k)$. However, it is not true that $[n]_q!=\sum_{k=0}^n c_q[n,k]$, where $c_q[n,k]=S_{1,q}[n,k]$. It would be interesting to derive an expression for $[n+m]_q!$ similar to Mez\"o's for $q$-Stirling numbers of the first kind.

We also ask if a modification of the rook model for the Type II generalized $q$-Stirling numbers can be made to given another combinatorial interpretation for both the Type I and Type II $p,q$-analogues. It is also possible that the techniques outlined in this paper may be modified to obtain the corresponding identities for the generalized Stirling numbers that arise from other rook models, such as those in \cite{RemWac} and \cite{Brig}.

Lastly, we note that other forms of convolution identities are given in \cite[Theorem 1]{AgoDil} and \cite[page 5]{MedLer}. Can these identities and their $S_{s,q}[n,k]$ versions be proved using partitions on rook boards?

\bibliographystyle{amsplain}

\newpage

\end{document}